\documentclass[oneside, 12pt]{amsart}

\usepackage{amsfonts}
\usepackage{amssymb}
\usepackage{latexsym}
\usepackage{graphics}
\usepackage[2emode]{psfrag}
\usepackage{mathrsfs}
\usepackage{amsthm}
\usepackage{amsmath}
\usepackage[all]{xy}
\usepackage{enumerate}
\usepackage{hyperref}

\begin{document}

\newcommand{\thlabel}[1]{\label{th:#1}}
\newcommand{\thref}[1]{Theorem~\ref{th:#1}}
\newcommand{\selabel}[1]{\label{se:#1}}
\newcommand{\seref}[1]{Section~\ref{se:#1}}
\newcommand{\lelabel}[1]{\label{le:#1}}
\newcommand{\leref}[1]{Lemma~\ref{le:#1}}
\newcommand{\prlabel}[1]{\label{pr:#1}}
\newcommand{\prref}[1]{Proposition~\ref{pr:#1}}
\newcommand{\colabel}[1]{\label{co:#1}}
\newcommand{\coref}[1]{Corollary~\ref{co:#1}}
\newcommand{\relabel}[1]{\label{re:#1}}
\newcommand{\reref}[1]{Remark~\ref{re:#1}}
\newcommand{\exlabel}[1]{\label{ex:#1}}
\newcommand{\exref}[1]{Example~\ref{ex:#1}}
\newcommand{\delabel}[1]{\label{de:#1}}
\newcommand{\deref}[1]{Definition~\ref{de:#1}}
\newcommand{\eqlabel}[1]{\label{eq:#1}}
\newcommand{\equref}[1]{(\ref{eq:#1})}

\newcommand{\Hom}{{\rm Hom}}
\newcommand{\End}{{\rm End}}
\newcommand{\Ext}{{\rm Ext}}
\newcommand{\Fun}{{\rm Fun}}
\newcommand{\Mor}{{\rm Mor}\,}
\newcommand{\Aut}{{\rm Aut}\,}
\newcommand{\Hopf}{{\rm Hopf}\,}
\newcommand{\Ann}{{\rm Ann}\,}
\newcommand{\Ker}{{\rm Ker}\,}
\newcommand{\Coker}{{\rm Coker}\,}
\newcommand{\im}{{\rm Im}\,}
\newcommand{\coim}{{\rm Coim}\,}
\newcommand{\Trace}{{\rm Trace}\,}
\newcommand{\Char}{{\rm Char}\,}
\newcommand{\Mod}{{\bf mod}}
\newcommand{\Spec}{{\rm Spec}\,}
\newcommand{\Span}{{\rm Span}\,}
\newcommand{\sgn}{{\rm sgn}\,}
\newcommand{\Id}{{\rm Id}\,}
\newcommand{\Com}{{\rm Com}\,}
\newcommand{\codim}{{\rm codim}}
\newcommand{\Mat}{{\rm Mat}}
\newcommand{\Coint}{{\rm Coint}}
\newcommand{\Incoint}{{\rm Incoint}}
\newcommand{\can}{{\sf can}}
\newcommand{\sign}{{\rm sign}}
\newcommand{\kar}{{\rm kar}}
\newcommand{\rad}{{\rm rad}}
\newcommand{\ev}{{\sf ev}}
\newcommand{\m}{{\sf m}}
\newcommand{\sd}{{\sf d}}
\def\colim{{\rm colim}\,}

\def\Ab{\underline{\underline{\rm Ab}}}
\def\lan{\langle}
\def\ran{\rangle}
\def\ot{\otimes}

\def\id{\textrm{{\small 1}\normalsize\!\!1}}
\def\To{{\multimap\!\to}}
\def\bigperp{{\LARGE\textrm{$\perp$}}} 
\newcommand{\QED}{\hspace{\stretch{1}}
\makebox[0mm][r]{$\Box$}\\}

\def\AA{{\mathbb A}}
\def\BB{{\mathbb B}}
\def\CC{{\mathbb C}}
\def\DD{{\mathbb D}}
\def\EE{{\mathbb E}}
\def\FF{{\mathbb F}}
\def\GG{{\mathbb G}}
\def\HH{{\mathbb H}}
\def\II{{\mathbb I}}
\def\JJ{{\mathbb J}}
\def\KK{{\mathbb K}}
\def\LL{{\mathbb L}}
\def\MM{{\mathbb M}}
\def\NN{{\mathbb N}}
\def\OO{{\mathbb O}}
\def\PP{{\mathbb P}}
\def\QQ{{\mathbb Q}}
\def\RR{{\mathbb R}}
\def\SS{{\mathbb S}}
\def\TT{{\mathbb T}}
\def\UU{{\mathbb U}}
\def\VV{{\mathbb V}}
\def\WW{{\mathbb W}}
\def\XX{{\mathbb X}}
\def\YY{{\mathbb Y}}
\def\ZZ{{\mathbb Z}}

\def\aa{{\mathfrak A}}
\def\bb{{\mathfrak B}}
\def\cc{{\mathfrak C}}
\def\dd{{\mathfrak D}}
\def\ee{{\mathfrak E}}
\def\ff{{\mathfrak F}}
\def\gg{{\mathfrak G}}
\def\hh{{\mathfrak H}}
\def\ii{{\mathfrak I}}
\def\jj{{\mathfrak J}}
\def\kk{{\mathfrak K}}
\def\ll{{\mathfrak L}}
\def\mm{{\mathfrak M}}
\def\nn{{\mathfrak N}}
\def\oo{{\mathfrak O}}
\def\pp{{\mathfrak P}}
\def\qq{{\mathfrak Q}}
\def\rr{{\mathfrak R}}
\def\ss{{\mathfrak S}}
\def\tt{{\mathfrak T}}
\def\uu{{\mathfrak U}}
\def\vv{{\mathfrak V}}
\def\ww{{\mathfrak W}}
\def\xx{{\mathfrak X}}
\def\yy{{\mathfrak Y}}
\def\zz{{\mathfrak Z}}

\def\aaa{{\mathfrak a}}
\def\bbb{{\mathfrak b}}
\def\ccc{{\mathfrak c}}
\def\ddd{{\mathfrak d}}
\def\eee{{\mathfrak e}}
\def\fff{{\mathfrak f}}
\def\ggg{{\mathfrak g}}
\def\hhh{{\mathfrak h}}
\def\iii{{\mathfrak i}}
\def\jjj{{\mathfrak j}}
\def\kkk{{\mathfrak k}}
\def\lll{{\mathfrak l}}
\def\mmm{{\mathfrak m}}
\def\nnn{{\mathfrak n}}
\def\ooo{{\mathfrak o}}
\def\ppp{{\mathfrak p}}
\def\qqq{{\mathfrak q}}
\def\rrr{{\mathfrak r}}
\def\sss{{\mathfrak s}}
\def\ttt{{\mathfrak t}}
\def\uuu{{\mathfrak u}}
\def\vvv{{\mathfrak v}}
\def\www{{\mathfrak w}}
\def\xxx{{\mathfrak x}}
\def\yyy{{\mathfrak y}}
\def\zzz{{\mathfrak z}}

\newcommand{\aA}{\mathscr{A}}
\newcommand{\bB}{\mathscr{B}}
\newcommand{\cC}{\mathscr{C}}
\newcommand{\dD}{\mathscr{D}}
\newcommand{\eE}{\mathscr{E}}
\newcommand{\fF}{\mathscr{F}}
\newcommand{\gG}{\mathscr{G}}
\newcommand{\hH}{\mathscr{H}}
\newcommand{\iI}{\mathscr{I}}
\newcommand{\jJ}{\mathscr{J}}
\newcommand{\kK}{\mathscr{K}}
\newcommand{\lL}{\mathscr{L}}
\newcommand{\mM}{\mathscr{M}}
\newcommand{\nN}{\mathscr{N}}
\newcommand{\oO}{\mathscr{O}}
\newcommand{\pP}{\mathscr{P}}
\newcommand{\qQ}{\mathscr{Q}}
\newcommand{\rR}{\mathscr{R}}
\newcommand{\sS}{\mathscr{S}}
\newcommand{\tT}{\mathscr{T}}
\newcommand{\uU}{\mathscr{U}}
\newcommand{\vV}{\mathscr{V}}
\newcommand{\wW}{\mathscr{W}}
\newcommand{\xX}{\mathscr{X}}
\newcommand{\yY}{\mathscr{Y}}
\newcommand{\zZ}{\mathscr{Z}}

\newcommand{\Aa}{\mathcal{A}}
\newcommand{\Bb}{\mathcal{B}}
\newcommand{\Cc}{\mathcal{C}}
\newcommand{\Dd}{\mathcal{D}}
\newcommand{\Ee}{\mathcal{E}}
\newcommand{\Ff}{\mathcal{F}}
\newcommand{\Gg}{\mathcal{G}}
\newcommand{\Hh}{\mathcal{H}}
\newcommand{\Ii}{\mathcal{I}}
\newcommand{\Jj}{\mathcal{J}}
\newcommand{\Kk}{\mathcal{K}}
\newcommand{\Ll}{\mathcal{L}}
\newcommand{\Mm}{\mathcal{M}}
\newcommand{\Nn}{\mathcal{N}}
\newcommand{\Oo}{\mathcal{O}}
\newcommand{\Pp}{\mathcal{P}}
\newcommand{\Qq}{\mathcal{Q}}
\newcommand{\Rr}{\mathcal{R}}
\newcommand{\Ss}{\mathcal{S}}
\newcommand{\Tt}{\mathcal{T}}
\newcommand{\Uu}{\mathcal{U}}
\newcommand{\Vv}{\mathcal{V}}
\newcommand{\Ww}{\mathcal{W}}
\newcommand{\Xx}{\mathcal{X}}
\newcommand{\Yy}{\mathcal{Y}}
\newcommand{\Zz}{\mathcal{Z}}

\def\units{{\mathbb G}_m}
\def\rightact{\hbox{$\leftharpoonup$}}
\def\leftact{\hbox{$\rightharpoonup$}}

\def\*C{{}^*\hspace*{-1pt}{\Cc}}
\def\*c{{}^*\hspace*{-1pt}{\cc}}

\def\text#1{{\rm {\rm #1}}}

\def\smashco{\mathrel>\joinrel\mathrel\triangleleft}
\def\cosmash{\mathrel\triangleright\joinrel\mathrel<}

\def\ol{\overline}
\def\ul{\underline}
\def\dul#1{\underline{\underline{#1}}}
\def\Nat{\dul{\rm Nat}}
\def\Set{\dul{\rm Set}}

\renewcommand{\subjclassname}{\textup{2000} Mathematics Subject
     Classification}

\newtheorem{proposition}{Proposition}[section] 
\newtheorem{lemma}[proposition]{Lemma}
\newtheorem{corollary}[proposition]{Corollary}
\newtheorem{theorem}[proposition]{Theorem}

\theoremstyle{definition}
\newtheorem{Definition}[proposition]{Definition}
\newtheorem{example}[proposition]{Example}
\newtheorem{examples}[proposition]{Examples}

\theoremstyle{remark}
\newtheorem{remarks}[proposition]{Remarks}
\newtheorem{remark}[proposition]{Remark}

\title[Equivalences between categories of modules and comodules]{Equivalences between categories of modules and categories of comodules}
 \author{Joost Vercruysse}    
 \address{Vrije Universiteit Brussel VUB, Pleinlaan 2, B-1050,
  Brussel, Belgium} 
 \email{joost.vercruysse@vub.ac.be}   
 \urladdr{http://homepages.vub.ac.be/\~{}jvercruy/}
 \date{August, 2006} 

\begin{abstract}
We show the close connection between appearingly different Galois theories for comodules introduced recently in \cite{GTV} and \cite{Wis:galcom}. Furthermore we study equivalences between categories of comodules over a coring and modules over a firm ring. We show that these equivalences are related to Galois theory for comodules.
\end{abstract}

\maketitle

\section{Introduction}\selabel{intro}

The classical Galois theory of for finite field extensions has a formulation in terms of Hopf algebras. If one takes $H=kG^*$ where $G$ is the group of automorphisms of a field $L\supset k$, and puts $F=L^G$ the field fixed by $G$, than $L$ is an $H$-comodule and $L$ is a Galois extension of $F$ if and only if $L\otimes_FL\cong L\otimes_kH$, by a canonical map (see e.g. \cite[Example 6.4.3 1)]{DNR}). 
This example has led to a range of generalizations, all under the name of Galois theory (We refer to \cite{Cae:desc} and \cite{Wis:galfield} for a profound overview). 
A comodule algebra $A$ over the Hopf algebra $H$, is called an $H$-Galois extension of 
$A^{{\rm co}H}$ if and only if the canonical map
$$\can:A\otimes_{A^{{\rm co}H}}A\to A\otimes_k H,\quad\can(a\otimes a')=aa'_{[0]}\otimes a'_{[1]},$$
is an isomorphism,
where $\rho^A(a)=a_{[0]}\ot a_{[1]}$ denotes the $H$-coaction on $A$.

It was observed in \cite{Brz:strcor} that Hopf-Galois theory and its generalizations can be beautifully reformulated using the language of corings. The Galois corings of \cite{Brz:strcor} have been generalized by El Kaoutit and G\'omez-Torrecillas in \cite{EGT:comcor}, leading to the notion of Galois comodule. An essential aspect of a Galois comodule from \cite{EGT:comcor} is that it is necessarily finite (by this, we mean that it is finitely generated and projective as a right module). 

In recent publications, several attempts have been made to introduce a notion of infinite Galois comodule. In \cite{GTV}, G\'omez-Torrecillas and the author consider Galois comodules over firm rings, this construction generalizes infinite Galois comodules introduced in \cite{EGT:infinite, CDV:colimit}. An alternative generalization has been proposed by Wisbauer in \cite{Wis:galcom}. The aim of this paper is to study the relationship between these two notions (\seref{comparision}). Our results are based on characterisations of equivalences between categories of comodules and categories of (firm) modules (\seref{equivalence}).

In \seref{firm} we have collected basic properties of Galois comodules over firm rings. One remarkable fact is the observation that firm rings and their firm modules can be identified with corings over the Dorroh extension of this firm ring and the comodules over this coring (\thref{firmcoring}). 
%
Based on properties of pairs of adjoint functors 
between categories of modules and comodules over firm rings, we introduce in \seref{equivalence} comonadic-Galois comodules as a refinement of the notion of Galois comodule in the sense of \cite{Wis:galcom}. It follows from elementary observations that comonadic Galois-comodules are 
Galois comodules in the sense of \cite{Wis:galcom}. Nevertheless, to study equivalences between comodule and module categories, the interaction with firm rings is essential to pass the finiteness border, as follows from the discussion below.
Our main result (\thref{equiv}) gives a characterisation of pairs of inverse equivalences between the category of modules over a firm ring and the category of comodules over a coring. More specifically, any such equivalence can be obtained from a Galois comodule in the sense of \cite{GTV}. As an application, we characterize equivalences between categories of comodules and modules over rings with unit, local units or idempotent local units (see Corollaries
\ref{co:locunits}, \ref{co:locunits} and \ref{co:equivunit}). 
An important implication of our theory is the fact that any equivalence between a category of comodules and a category of modules implies a certain generalised projectivity condition (what we named `firmly projectivty' in this note) on a generator of the category of comodules, which plays the role of the Galois comodule. In case of an equivalence between a category of comodules with a category of modules over a ring with unit, this implies that the generator is always finitely generated and projective, and `finite' Galois theory comes into play. This indicates that firm rings play an important role if one wants to study Galois theories that go beyond the case where the Galois comodule is finitely generated and projective.
The main result of \seref{comparision} states that the notion of Galois comodule and comonadic Galois comodule coincide when the comodule under considaration satisfies a generalized notion of projectivity (see \thref{comparmain}). Another remarkable result is \thref{compar2} : if $\Sigma$ is a comonadic Galois comodule in the sense of \cite{Wis:galcom}, then the underlying coring $\cc$ is isomorphic to a comatrix coring associated to a comatrix coring context over a firm ring, as introduced in \cite{GTV}.

\section{Galois theory for comodules bounded by firm rings}\selabel{firm}

\subsection*{Corings and Comodules}
Let $A$ be an associative ring (with unit). An $A$-\emph{coring} consists of an $A$-bimodule $\cc$, and two $A$-bilinear maps
$$\Delta_\cc:\cc\to\cc\otimes_A\cc, \quad \varepsilon_\cc :\cc\to A,$$
such that the following diagrams commute
\[
\xymatrix{
\cc \ar[rr]^-{\Delta_\cc} \ar[d]_{\Delta_\cc} && \cc\otimes_A\cc \ar[d]^{\Delta_\cc\otimes_A\cc}\\
\cc\otimes_A\cc \ar[rr]_-{\cc\otimes_A\Delta_\cc} && \cc\otimes_A\cc\otimes_A\cc
}
\qquad
\xymatrix{
\cc \ar[rrd]^-{\Delta_\cc} && \cc\otimes_AA \ar[ll]_-{\cong}\\
A\otimes_A\cc \ar[u]^{\cong} && \cc\otimes_A\cc \ar[ll]^{\varepsilon_\cc\otimes_A\cc} \ar[u]_{\cc\otimes_A\varepsilon_\cc} 
}
\]
Remark our convention to denote the identity morphism on an object $X$ by the same character $X$.
We will make use of the Sweedler notation $\Delta(c)=c_{(1)}\otimes_Ac_{(2)}$ and $(\cc\otimes_A\Delta_\cc)\circ\Delta(c)
=(\Delta_\cc\otimes_A\cc)\circ\Delta(c)=c_{(1)}\otimes_Ac_{(2)}\otimes_Ac_{(3)}$.
 
A right $\cc$-comodule is a pair $(M,\rho_M)$, where $M$ is a right $A$-module and $\rho:M\to M\ot_A\cc, \rho_M(m)=m_{[0]}\ot_Am_{[1]}$ is a right $A$-linear map, such that $(\rho_M\ot_A\cc)\circ\rho_M=(M\ot_A\Delta_\cc)\circ\rho_M$ and $(M\ot_A\varepsilon_\cc)\circ\rho_M=M$. A right $A$-linear map $f:M\to N$ between the right $\cc$-comodules $M$ and $N$ is said to be right $\cc$-colinear if $\rho_N\circ f=(f\ot_A\cc)\circ\rho_M$. The category of all right $\cc$-comodules and right $\cc$-colinear maps is denoted by $\Mm^\cc$. If $B$ is a second ring, then ${_B\Mm^\cc}$ is the notation for the category of all left $B$-modules $M$ which also have a right $\cc$-comodule structure and such that the comultiplication $\rho_M$ on $M$ is left $B$-linear. Morphisms in ${_B\Mm^\cc}$ are maps that are left $B$- and right $\cc$-colinear. 

For a comprehensive introduction to the theory of corings an comodules we refer to the monograph \cite{BrzWis:book}.

\subsection*{Firm rings}

Let $R$ be an arbitrary ring, not necessary with unit. By $\Mm_R$ we denote the category of right $R$-linear maps and all right $R$-modules satisfying the property that the multiplication map
$$\mu_{M,R}:M\otimes_RR\to M,\qquad \mu_{M,R}(m\otimes_Rr)=mr,$$
establishes an isomorphism. In this case, we denote the inverse map as
$${\sf d}_{M,R}:M\to M\otimes_RR,\qquad {\sf d}_{M,R}(m)=m^r\otimes r.$$
We call $\Mm_R$ the category of \emph{firm right $R$-modules}. Similary one introduces the category ${_R\Mm}$. It is easy to check that $R\in\Mm_R$ is equivalent to $R\in{_R\Mm}$, in this situation we refer to $R$ as a \emph{firm ring} and we will denote the structure maps of $R$ by $\mu_R$ and $\sd_R$. The category $\Mm_R$ is abelian provided that $R\in{\Mm_R}$, in that case $R$ is also a generator of the category. We call $M\in{\Mm_R}$ flat provided that the functor $M\ot_R-:\widetilde{\Mm}_R\to\Ab$ is exact, where $\widetilde{\Mm}_R$ denotes the category of all (possibly non-firm) $R$-modules. 
If $R$ is flat as a left $R$-module, then $\Mm_R$ is a Grothendieck category and kernels, cokernels and coproducts can be computed in $\Ab$. 

Let $R$ be a (non-unital) ring. Recall (see e.g. \cite[section 1.5]{Wis:book}) that the Dorroh-extension of $R$ is a unital ring $\widehat{R}=R\times\ZZ$ containing $R$ as a two-sided ideal by the canonical injection $\iota:R\to \widehat{R},\ \iota(r)=(r,0)$. The multiplication in 
$\widehat{R}$ is given by
$$(r,x)(r',x')=(rr'+rx'+xr',xx')$$
for all $r,r'\in R$ and $x,x'\in \ZZ$ and $(0,1)$ is the unit of $\widehat{R}$. 
For any $M\in\widetilde{\Mm}_R$, we can define a unital right $\widehat{R}$-action via
$$m(r,x)=mr+mx,$$
for all $m\in M$, $r\in R$ and $x\in\ZZ$. In this way, we obtain an isomorphism of categories between between $\widetilde{\Mm}_R$ and $\Mm_{\widehat{R}}$. 
In general, there exists no unital ring $S$ such that $\Mm_R$ is isomorphic to $\Mm_S$, even if $R$ is firm. However, we can describe $\Mm_R$ as a category of comodules over a coring (with a unital base ring).
Remark first that the canonical surjection $M\ot_RR\to M\ot_{\widehat{R}}R$ is an isomorphism.

\begin{lemma}\lelabel{2.1}
Let $R$ be a non-unital ring. Then $R$ is a firm ring if and only if there exists a bimodule map 
$\Delta:R\to R\ot_RR\cong R\ot_{\widehat{R}}R$ such that $(R,\Delta,\iota)$ is an $\widehat{R}$-coring.
\end{lemma}

\begin{proof}
Suppose $R$ is a firm ring. Take $\Delta=\sd_R$. The coassociativity of $\Delta$ follows by the associativity of $\mu_R$, since they are two-sided inverses in $\Mm_R$. Denote $\Delta(r)=r_{(1)}\ot_{\widehat{R}}r_{(2)}$, then $\iota(r_{(1)})r_{(2)}=r_{(1)}\iota(r_{(2)})=r_{(1)}r_{(2)}=\mu_R\circ\sd_R(r)=r$, so the counit property holds as well.

Conversely, if $R$ is an $\widehat{R}$-coring, one can easily check 
\begin{eqnarray*}
(\mu_R\circ\Delta)(r)&=&r_{(1)}r_{(2)}=r_{(1)}\iota(r_{(2)})=r\\
(\Delta\circ\mu_R)(\sum_i r_i\ot_Rr'_i)&=&\Delta(\sum_ir_ir'_i)=\sum_i\Delta(r_i)r'_i\\
&=&\sum_ir_{i(1)}\ot_Rr_{i(2)}r'_i=\sum_ir_{i(1)}r_{i(2)}\ot_Rr'_i\\
&=&\sum_ir_i\ot_Rr'_i
\end{eqnarray*}
so $\Delta$ is a two-sided inverse for $\mu_R$, so $R$ is a firm ring.
\end{proof}

\begin{theorem}\thlabel{firmcoring}
Assume that $R$ is a firm ring and let $\rr=(R,\Delta,\iota)$ be the associated $\widehat{R}$-coring from \leref{2.1}. Then the categories $\Mm_R$ and $\Mm^\rr$ are isomorphic.
\end{theorem}

\begin{proof}
Take $(M,\rho)\in\Mm^\rr$. Then $M$ is in particular a right $\widehat{R}$-module, so a non-unital $R$-module. Let $\sd_{M,R}$ be the composition $M\to M\ot_{\widehat{R}}R\cong M\ot_RR$. As in \leref{2.1} we can easily check that $\mu_{M,R}$ and $\sd_{M,R}$ are mutual inverses, so $M$ is firm as a right $R$-module.

Conversely, let $M$ be a firm right $R$-module and define $\rho$ as the composition
\[
\xymatrix{
M\ar[r]^-{\sd_{M,R}} & M\ot_RR \ar[r]^\cong & M\ot_{\widehat{R}}R
}
\]
The coassociativity of $\rho$ follows by the coassociativity of $\mu_{M,R}$. Finally, for any $m\in M$,
$$(M\ot_R\iota)\circ\rho(m)= m^r\iota(r)=m^rr=m,$$
so the counit property is satisfied as well.
\end{proof}

\subsection*{Firmly projective modules}

Let $A$ and $B$ be rings, $\Sigma\in{_B\Mm_A}$, $\Sigma^\dagger\in{_A\Mm_B}$ and $\mu:\Sigma^\dagger\ot_B\Sigma\to A$ an $A$-bilinear map. Then we call $(\Sigma,\Sigma^\dagger,\mu)$ a \emph{dual pair}. Recall the construction of an elementary $B$-ring $Z^\dagger=\Sigma\ot_A\Sigma^\dagger$ associated to such a dual pair. The multiplication of this ring is given by $$\mu_{Z^\dagger}=\Sigma\ot_A\mu\ot_A\Sigma^\dagger.$$ In general $Z^\dagger$ has no unit. $\Sigma$ and $\Sigma^\dagger$ are a left, respectively right $Z^\dagger$-module with actions given by
$$\mu_{Z^\dagger,\Sigma}=\Sigma\ot_A\mu,$$ 
and
$$\mu_{\Sigma^\dagger,Z^\dagger}=\mu\ot_A\Sigma^\dagger.$$

$(\Sigma,\Sigma^\dagger,\mu)$ is said to be an \emph{$R$-firm dual pair} if their exists a firm ring $R$ together with a ring morphism $\iota:R\to\Sigma\ot_A\Sigma^\dagger$, $\iota(r)=e_r\ot_Af_r$ (summation understood) such that $\Sigma$ and $\Sigma^\dagger$ are a firm left, respectively a firm right $R$-module, where the $R$-action is induced by $\iota$.

In particular, we can consider the notion of a \emph{firm dual pair}, taking $R=Z^\dagger$ and $\iota$ the identity.

Consider the dual pair $(\Sigma,\Sigma^*,\ev)$, where $\ev:\Sigma^*\ot_B\Sigma\to A$ is the evaluation map. We will denote the associated $B$-ring by $Z$, the multiplication in $Z$ as $\mu_Z$, and the actions of $Z$ on $\Sigma$ and $\Sigma^\dagger$ by $\mu_{Z,\Sigma}$ and $\mu_{\Sigma^\dagger,Z}$.

In this paper, we will call a $B$-$A$-bimodule $\Sigma$ \emph{firmly projective} if and only if the elementary $B$-ring $Z=\Sigma\ot_A\Sigma^*$ is a firm ring and $\Sigma$ is firm as a left $Z$-module.

Similary, for a firm $B$-ring $R$, $\Sigma$ will be named \emph{$R$-firmly projective} if and only if there exists a ring morphism $\iota : R\to Z$, $\iota(r)=e_r\ot_Af_r$ and $\Sigma$ is a firm $R$-module under the $R$-module structure induced by $\iota$.

Obviously, if $\Sigma$ is firmly projective, then it is $Z$-firmly projective. 
In the following Proposition we have collected elementary results about $R$-firmly projective modules and $R$-dual pairs.

\begin{proposition}\lelabel{firm}
\begin{enumerate}[(i)]
\item If $\Sigma$ is $R$-firmly projective for any firm $B$-ring $R$, then $\Sigma$ is firmly projective.
\item If $\Sigma$ is $R$-firmly projective, then $(\Sigma,\Sigma^\dagger,\mu)$ is a $R$-firm dual pair, with $\Sigma^\dagger=\Sigma^*\ot_RR$ and $\mu=\ev\circ(\Sigma^*\ot_R\mu_{Z,\Sigma})$.
\item If $(\Sigma,\Sigma^\dagger,\mu)$ is an $R$-firm dual pair, then $\Sigma$ is $R$-firmly projective and $\Sigma^\dagger\cong\Sigma^*\ot_RR$ as $A$-$R$-bimodule.
\end{enumerate}
\end{proposition}

\begin{proof}
\ul{(i)}. First remark that the multiplication $R\ot_B\Sigma\to\Sigma$ is a right $A$-linear (and thus an $R$-$A$ bilinear) map, as it is the composition of right $A$-linear maps: $(\Sigma\ot_A\ev)\circ(\iota\ot_B\Sigma)$. This implies that if $\Sigma$ is $R$-firmly projective and, a fortiori, firm as a left $R$-module, the isomorphism $R\ot_R\Sigma\cong\Sigma$ holds as an isomorphism of $R$-$A$ bimodules. This isomorphism induces a well-defined map
\begin{equation*}
\sd:\Sigma\ot_A\Sigma^* \cong R\ot_R\Sigma\ot_A\Sigma^* \stackrel{\iota\ot_RZ}{\longrightarrow} \Sigma\ot_A\Sigma^*\ot_R\Sigma\ot_A\Sigma^* \stackrel{\pi}{\longrightarrow} \Sigma\ot_A\Sigma^*\ot_Z\Sigma\ot_A\Sigma^*,
\end{equation*}
or $\sd(u\ot_A\varphi)=e_r\ot_Af_r\ot_Zu^r\ot_A\varphi$ for all $u\ot_A\varphi\in\Sigma\ot_A\Sigma^*$. Obviously, $\mu_Z\circ\sd=Z$, since $e_rf_r(u^r)\ot_A\varphi=u\ot_A\varphi$. Conversely, consider an element $u\ot_A\varphi\ot_Zv\ot_A\psi\in\Sigma\ot_A\Sigma^*\ot_Z\Sigma\ot_A\Sigma^*$ (summation understood). Then we compute that 
\begin{eqnarray*}
(\sd\circ\mu_Z)(u\ot_A\varphi\ot_Zv\ot_A\psi)&=&
e_r\ot_Af_r\ot_Z(u\varphi(v))^r\ot_A\psi\\
&=&e_r\ot_Af_r\ot_Zu^r\varphi(v)\ot_A\psi\\
&=&e_r\ot_Af_r\ot_Z(u^r\ot_A\varphi)\cdot v\ot_A\psi\\
&=&e_r\ot_Af_r\cdot(u^r\ot_A\varphi)\ot_Zv\ot_A\psi\\
&=&e_r\ot_Af_r(u^r)\varphi\ot_Zv\ot_A\psi\\
&=&e_rf_r(u^r)\ot_A\varphi\ot_Zv\ot_A\psi\\
&=&u\ot_A\varphi\ot_Zv\ot_A\psi,
\end{eqnarray*}
so $\sd\circ\mu_Z=Z$ and $Z$ is a firm ring. In a similar way, $\Sigma$ is a firm $Z$-module.

\ul{(ii)} follows immediately since $M\ot_RR$ is a firm right $R$-module for any right $R$-module $M$.

\ul{(iii)}. Consider the morphism $\zeta:\Sigma^\dagger\to\Sigma^*$, $\zeta(\varphi)(u)=\mu(\varphi\ot_Bu)$ for any $\varphi\in\Sigma^\dagger$ and $u\in\Sigma$. Recall that we have a morphism $\iota:R\to \Sigma\ot_A\Sigma^\dagger, \iota(r)=e_r\ot_Rf_r$.
Then we can consider the morphisms
\begin{eqnarray*}
\alpha:\Sigma^\dagger \to \Sigma^*\ot_RR,&& \alpha(\varphi)=\zeta(\varphi^r)\ot_Rr\\
\beta:\Sigma^*\ot_RR\to \Sigma^\dagger,&& \beta(\psi\ot_Rr)=\psi(e_r)f_r
\end{eqnarray*}
A similar calculation as in part $(i)$ shows that $\alpha$ and $\beta$ are each other inverses.
\end{proof}

Examples of firmly projective modules and firm dual pairs can be easily obtained from (locally) projective modules.

\subsection*{Galois comodules}

Recall from \cite{BGT} the notion of a \emph{comatrix coring context} (there for rings with unit). Let $A$ and $R$ be firm rings, $\Sigma\in{_R\Mm_A}$, $\Sigma^\dagger\in{_A\Mm_R}$ and consider two bilinear maps $\eta:R\to\Sigma\ot_A\Sigma^\dagger$ and $\varepsilon:\Sigma^\dagger\ot_R\Sigma\to A$. Then $(R,A,\Sigma,\Sigma^\dagger,\eta,\varepsilon)$ is a comatrix coring context if and only if the following diagrams commute.
\begin{equation}\eqlabel{ccc}
\xymatrix{
\Sigma \ar[rr]^-{\cong} \ar[d]_{\cong} && R\ot_R\Sigma \ar[d]^{\eta\ot_R\Sigma}\\
\Sigma\ot_AA && \Sigma\ot_A\Sigma^\dagger\ot_R\Sigma \ar[ll]^{\Sigma\ot_A\varepsilon}
}\qquad
\xymatrix{
\Sigma^\dagger \ar[rr]^-{\cong} \ar[d]_{\cong} && \Sigma^\dagger\ot_RR \ar[d]^{\Sigma^\dagger\ot_R\eta}\\
A\ot_A\Sigma^\dagger && \Sigma^\dagger\ot_R\Sigma\ot_A\Sigma^\dagger \ar[ll]^{\varepsilon\ot_A\Sigma^\dagger}
}
\end{equation}

In this situation, one can construct a $R$-ring $Z^\dagger=\Sigma\ot_A\Sigma^\dagger$ with unit $\eta$ and multiplication $\Sigma\ot_A\varepsilon\ot_A\Sigma^\dagger$ (remark that $(\Sigma,\Sigma^\dagger,\varepsilon)$ is a dual pair) and an $A$-coring $\dd=\Sigma^\dagger\ot_R\Sigma$ with counit $\varepsilon$ and comultiplication $\Sigma^\dagger\ot_B\eta\ot_B\Sigma$ (see \cite[Poposition 1.1]{CDV:colimit}). Moreover $\Sigma$ is a left $Z^\dagger$-module and a right $\dd$-comodule. $\Sigma^\dagger$ is a right $Z^\dagger$-module and a left $\dd$-comodule. Actions and coactions are given by the following formula
\begin{eqnarray*}
\rho_{\Sigma,\dd}=(\eta\ot_B\Sigma)\circ\sd_{B,\Sigma} && \mu_{Z^\dagger,\Sigma}=\mu_{\Sigma,A}\circ(\Sigma\ot_A\varepsilon)\\
\rho_{\dd,\Sigma^\dagger}=(\Sigma^\dagger\ot_B\eta)\circ\sd_{\Sigma^\dagger,B} && \mu_{\Sigma^\dagger,Z^\dagger}=\mu_{A,\Sigma^\dagger}\circ(\varepsilon\ot_A\Sigma^\dagger)
\end{eqnarray*}

\begin{theorem}\thlabel{firmadjoint}
Let $A$ and $R$ be firm rings, $\Sigma\in{_A\Mm_R}$ and $\Sigma^\dagger\in{_R\Mm_A}$. The following statements are equivalent
\begin{enumerate}[(i)]
\item $\Sigma$ is $R$-firmly projective;
\item there exists a comatrix coring context $(R,A,\Sigma,\Sigma^\dagger,\eta,\varepsilon)$;
\item there is a pair of adjoint functors $(-\ot_R\Sigma,-\ot_A\Sigma^\dagger)$ between the categories $\Mm_R$ and $\Mm_A$;
\item there exists a ``triple'' of adjoint funtors $(F,G,H)$ between the categories $\Mm_R$ and $\Mm_A$, i.e. $(F,G)$ is a pair of adjoint functors between $\Mm_R$ and $\Mm_A$ and $(G,H)$ is a pair of adjoint functors between $\Mm_A$ and $\Mm_R$, such that $\Sigma=F(R)$.
\end{enumerate}
\end{theorem}

\begin{proof}
$\ul{(i)\Rightarrow(ii)}$. Suppose $\Sigma$ is $R$-firmly projective for a firm $B$-ring $R$. Then one easily checks that $(R,A,\Sigma,\Sigma^\dagger,\iota,\mu)$ is a comarix coring context, where we denote as before, $\iota: R\to \Sigma\ot_A\Sigma^*$ the ring morphism from the definition of $R$-firmly projectivity, $\Sigma^\dagger=\Sigma^*\ot_RR$ and $\mu=\ev\circ(\Sigma^*\ot_R\mu_{Z,\Sigma})$.

$\ul{(ii)\Rightarrow(i)}$. The definition of a comatrix coring context implies that $\Sigma$ is a firm $R$-module. If we denote $\eta(r)=e_r\ot_Af_r$ for $r\in R$, then the commutativity of the left diagram of \equref{ccc} means that $r\cdot u=e_rf_r(u)$, which means exactly that the $R$-module structure of $\Sigma$ can be seen as induced by $\eta$, exactly as in the definition of an $R$-firmly projective module. 

$\ul{(ii)\Rightarrow(iii)}$. We restrict ourself to give the unit and counit of the adjunction and leave the other verifications to the reader
\begin{eqnarray*}
\alpha_N:N\to N\ot_R\Sigma\ot_A\Sigma^\dagger;&&\alpha(n)=n^r\ot_Re_r\ot_Af_r;\\
\beta_M:M\ot_A\Sigma^\dagger\ot_R\Sigma\to M;&&\beta(m\ot_A\varphi\ot_Ru)=m\varphi(u);
\end{eqnarray*}
for any $N\in\Mm_R$ and $M\in\Mm_A$.

$\ul{(iii)\Rightarrow (ii)}.$
Denote by $\alpha$ and $\beta$ the unit and counit of the adjunction. 
The unit evaluated in $R$ induces a right $R$-linear map
$$\alpha_R:R\to \Sigma\ot_A\Sigma^\dagger,$$
which is also left $R$-linear by naturality. 
Moreover, 
the counit evaluated in $A$ provides us with a map 
$$\beta_A:\Sigma^\dagger\ot_R\Sigma\to A,$$
this map becomes again $A$-bilinear by naturality.
We prove first that for any $N\in\Mm_R$, 
\begin{equation}\eqlabel{alphanatural}
\alpha_N\cong N\ot_R\alpha_R,
\end{equation}
i.e. the following diagram commutes.
\[
\xymatrix{
N \ar[rr]^-{\alpha_N} \ar[d]_\cong && N\ot_R\Sigma\ot_A\Sigma^\dagger \\
N\ot_RR \ar[rru]_-{N\ot_R\alpha_R}
}
\]
Take any $n\in N$, and consider the following map in $\Mm_R$, $f_n:R\to N$, $f_n(r)=nr$. By naturality of $\alpha$, the following diagram commutes,
\[
\xymatrix{
R \ar[rr]^-{\alpha_R} \ar[d]_{f_n} && R\ot_R\Sigma\ot_A\Sigma^\dagger \ar[d]^{f_n\ot_R\Sigma\ot_A\Sigma^\dagger}\\
N \ar[rr]^-{\alpha_N} && N\ot_R\Sigma\ot_A\Sigma^\dagger
}
\]
and \equref{alphanatural} follows.
Similarly we prove that for $M\in\Mm_A$, $\beta_M\cong M\ot_A\beta_A$.
The condition that $(F,G)$ is an adjoint pair means exactly that $\beta_{FM}\circ F(\alpha_M)=F(M)$ and $G(\beta_M)\circ\alpha_{GM}=G(M)$. If we evaluate these conditions in $\Sigma$ and $\Sigma^\dagger$, we obtain the following commutative diagrams, 
\[
\xymatrix{
\Sigma\cong R\ot_R\Sigma \ar[rr]^-{\alpha_\Sigma\cong \alpha_R\ot_R\Sigma} \ar@{=}[rrd] && \Sigma\ot_A\Sigma^\dagger\ot_R\Sigma \ar[d]^-{\Sigma\ot_A\beta_A} \\
&& \Sigma \cong \Sigma\ot_AA
}\qquad
\xymatrix{
\Sigma^\dagger\cong A\ot_A\Sigma^\dagger \ar@{=}[rrd] && \Sigma^\dagger\ot_R\Sigma\ot_A\Sigma^\dagger \ar[ll]_-{\beta_{\Sigma^\dagger}\cong\beta_A\ot_A\Sigma^\dagger} \\
&& \Sigma^\dagger \cong \Sigma^\dagger\ot_RR \ar[u]_{\Sigma^\dagger\ot_R\alpha_R}
}
\]
This means exactly that $(R,A,\Sigma,\Sigma^\dagger,\alpha_R,\beta_A)$ is a comatrix coring context.

$\ul{(iii)\Leftrightarrow(iv)}$. This follows immediately from the Eilenberg-Watts theorem (see also \thref{eilenbergwattsfirm})
\end{proof}

Let $(R,A,\Sigma,\Sigma^\dagger,\eta,\varepsilon)$ be a comatrix coring context as in \thref{firmadjoint}. Then $\Sigma$ is $R$-firmly projective and $\Sigma^\dagger=\Sigma^*\ot_RR$.
For the associated comatrix coring the following isomorphism holds $\dd=\Sigma^\dagger\ot_R\Sigma\cong\Sigma^*\otimes_R\Sigma$.
Since this will be useful in the sequel, we give the explicit formulas for the coactions of the associated coring and comodules.
\begin{eqnarray*}
\varepsilon_\dd : \Sigma^*\otimes_R\Sigma\to A,&& \varepsilon_\dd(\varphi \otimes_Ru)=\varphi(u)\\
\Delta_\dd : \Sigma^*\otimes_R\Sigma\to \Sigma^*\otimes_R\Sigma\otimes_A\Sigma^*\otimes_R\Sigma,&& 
\Delta_\dd(\varphi\otimes_Ru)=\varphi\otimes_R e_r\otimes_A f_r \otimes_R u^r\\
\rho_\Sigma : \Sigma \to \Sigma\otimes_A\Sigma^*\otimes_R\Sigma,& & \rho_\Sigma(u)=e_r\otimes_Af_r\otimes_Ru^r\\
\rho_{\Sigma^\dagger}:\Sigma^\dagger \to \Sigma^*\otimes_R\Sigma\otimes_A\Sigma^\dagger,&& \rho_{\Sigma^\dagger}(\varphi\otimes_Rr)=\varphi\otimes_Re_s\otimes_Af_s\otimes_Rr^s\nonumber,
\end{eqnarray*}
for any $u\in\Sigma$ and $\varphi\otimes_Rr\in\Sigma^\dagger=\Sigma\otimes_RR$ (summation understood).\\

Let $\cc$ be another $A$ coring and suppose that $\Sigma\in{_B\Mm_A}$. We will say that $\Sigma$ is $R$-firmly projective as $\cc$-comodule for a firm $B$-ring $R$ if there exists a ring morphism $\iota:R\to\Sigma\ot_A\Sigma^*$ and $\Sigma\in{_R\Mm^\cc}$, i.e. $\Sigma$ is a firm left $R$-module under the action induced by $\iota$ and the action of $R$ is $\cc$-colinear.

Under these conditions $\Sigma$ is called an $R\hbox{-}\cc$ Galois comodule if and only if the coring morphism
\begin{equation}\eqlabel{can}
\can:\Sigma^*\otimes_R\Sigma\to \cc,\quad \can(f\otimes_R u)=f(u_{[0]})u_{[1]}
\end{equation}
is an isomorphism. This definition of Galois comodule and its associated comatrix coring was given in \cite{GTV}.

\section{Equivalences beween categories of modules and comodules}\selabel{equivalence}

The classical Eilenberg-Watts theorem for adjoint functors between categories of modules over unital rings can be easily generalized to firm modules over firm rings. 

\begin{theorem}\thlabel{eilenbergwattsfirm}
Let $A$ and $B$ be two firm rings. For a functor $F:\Mm_B\to\Mm_A$, the following assertions are equivalent:
\begin{enumerate}[(i)]
\item $F$ has a right adjoint $G:\Mm_A\to\Mm_B$;
\item $F$ is right exact and preserves direct sums;
\item $F\cong - \ot_BP$ for some $B$-$A$ bimodule $P$.
\end{enumerate}
Moreover, when $(i)$-$(iii)$ hold, then the right adjoint $G$ is given by
$$G=\Hom_A(P,-)\ot_BB,$$
where $P=F(B)$ and $G$ is unique up to isomorphism.
\end{theorem}

\begin{proof}
The proof is identical to the classical proof (see e.g. \cite[Proposition 10.1]{Ste:book}). 
\end{proof}

Let $(F,G)$ be a pair of adjoint functors beween two module categories,
\[
\xymatrix{
\Mm_R \ar@<0.5ex>[rr]^-F && \Mm_A \ar@<0.5ex>[ll]^-G
}
\]
with unit $\eta$ and counit $\varepsilon$.
As it is well known, the composite functor $C=FG:\Mm_A\to \Mm_A$ defines a comonad (cotriple) $(C,F\eta G,\varepsilon)$ on the category $\Mm_A$ (see e.g. \cite[Chapter VI]{McLane}). Remark that an adjoint pair can be seen as a comatrix coring context in the bicategory of categories, functors and natural transformations. The construction of a comonad $C$ is then equivalent to the construction of a comatrix coring from this context within this bicategory. We can consider the subcategory $\Mm^C$ of $\Mm_A$ consisting of all `$C$-coalgebras'.
Recall that there exists a pair of adjoint functors
\[
\xymatrix{
\Mm_A \ar@<0.5ex>[rr]^-{G^C} && \Mm^C \ar@<0.5ex>[ll]^-{F^C}
}
\]
where $F^C$ is the forgetful functor and $G^C$ is given by the coristriction of the comonad functor $C$.
The original and new obtained adjoint pair of functors can be compared by a unique functor, such that we obain the following diagram
\begin{equation}\eqlabel{adjointcomonad}
\xymatrix{
\Mm_R \ar@<0.5ex>[rr]^-F \ar[rrdd]_-K && \Mm_A \ar@<0.5ex>[ll]^-G \ar@<0.5ex>[dd]^-{G^C}\\
\\
&& \Mm^C \ar@<0.5ex>[uu]^-{F^C}
}
\end{equation}
It turns out that the functor $K$ can be choosen as the corestriction of the functor $F$, since for any $M\in\Mm_R$, $F(M)\in\Mm^C$. 

If we apply the Eilenberg-Watts theorem to this situation, we obtain $$C\cong\Hom_A(\Sigma,-)\ot_R\Sigma,$$ 
Where we denote $\Sigma=F(R)$. By applying the Eilenberg-Watts theorem a second time, the comonad $C$ is of the form $-\ot_A\cc$ (i.e. induced by an $A$-coring $\cc$) if and only if the composite functor $C$ is right exact and preserves direct sums. If this is the case, then we find
\begin{equation}\eqlabel{comonadgalois}
-\ot_A\cc\cong\Hom_A(\Sigma,-)\ot_R\Sigma.
\end{equation}
By construction we find that $\Sigma\in{_R\Mm^\cc}$, i.e. there is a ring morphism $\jmath:R\to \End^\cc(\Sigma)$, and a natural transformation \equref{comonadgalois} takes the form
\begin{equation}\eqlabel{nattfgalois}
\can_M:\Hom_A(\Sigma,M)\ot_R\Sigma\to M\ot_A\cc,\qquad \can_M(f\ot_R u)=f(u_{[0]})u_{[1]}
\end{equation}
for all $M\in\Mm_A$.

For any firm ring $R$ and $A$-coring $\cc$, we will call an $R$-$\cc$ bicomodule $\Sigma$ an $R\hbox{-}\cc$ \emph{comonadic-Galois} comodule if \equref{comonadgalois} holds in the sense that \equref{nattfgalois} is an isomorphism for all $M\in\Mm_A$. Under the extra assumption that $\jmath$ is an isomorphism, this coincides with the notion of Galois comodule studied by Wisbauer in \cite{Wis:galcom}. If $\Sigma$ is an $R\hbox{-}\cc$ comonadic-Galois comodule, 
the diagram \equref{adjointcomonad} can be completed as follows
\[
\xymatrix{
\Mm_R \ar@<0.5ex>[rrrr]^-{-\ot_R\Sigma} \ar@<0.5ex>[rrrrdd]^-{-\ot_R\Sigma} &&&& \Mm_A \ar@<0.5ex>[llll]^-{\Hom_A(\Sigma,-)\ot_RR} \ar@<0.5ex>[dd]^-{-\ot_A\cc}\\
\\
&&&& \Mm^\cc \ar@<0.5ex>[uu]^-{F^\cc} \ar@<0.5ex>[lllluu]^-{\Hom^\cc(\Sigma,-)\ot_RR}
}
\]
The diagonal in this diagram is again an adjoint pair, this was proven in \cite[Lemma 4.2]{GTV}.
Let us prove a version of the Eilenberg-Watts theorem that states that any adjunction between a category of modules and comodules is of this form.

\begin{theorem}\thlabel{Eilenberg-Wattscoring}
Let $R$ be a firm ring and $\cc$ an $A$-coring. Let $F:\Mm_R\to\Mm^\cc$ be a functor and $\Sigma=F(R)$. Then the following assertions are equivalent:
\begin{enumerate}[(i)]
\item $F$ has a right adjoint $G:\Mm^\cc\to\Mm_R$;
\item $F$ is right exact and preserves direct sums;
\item $F\cong - \ot_R\Sigma$ for some $\Sigma\in{_R\Mm^\cc}$.\\ \vspace{-0.5cm}

\noindent\hspace{-1cm}
In this situation we also have the following properties
\item The right adjoint $G$ of $F$ is unique up to isomorphism and given by
\begin{equation}\eqlabel{expressionG}
G=\Hom^\cc(\Sigma,-)\ot_RR;
\end{equation}
\item there exists a ring morphism $\jmath:R\to \End^\cc(\Sigma)$;
\item there exists a natural transformation
\begin{equation}
\can_M:\Hom_A(\Sigma,M)\ot_R\Sigma\to M\ot_A\cc,
\end{equation}
given by $\can_M(f\ot_R u)=f(u_{[0]})u_{[1]}$.
\end{enumerate}
\end{theorem}

\begin{proof}
The proof goes allong the same lines as the classical Eilenberg-Watts theorem, but let us give a complete proof for sake of completeness.

$\ul{(i)\Rightarrow(ii)}$ is classical, $\ul{(iii)\Rightarrow(i)}$ and $(iv)$ follow from \cite[Lemma 4.2]{GTV}. 
$\ul{(ii)\Rightarrow(iii)}$. Put $\Sigma=F(R)$, then $\Sigma\in{_R\Mm^\cc}$. Since $R$ is a generator in $\Mm_R$, for any $M\in\Mm_R$ there exists an exact sequence,
\[
\xymatrix{
R^{(I)}\ar[rr]^-{f}&&R^{(J)}\ar[rr]^-{g}&&M\ar[rr]&&0.
}
\]
Here $(M,g)$ is exaclty the cokernel of $f$. 
Apply the functors $F$ and $-\ot_R\Sigma$ on this sequence. Since both functors are right exact and they  
preserve direct sums and cokernels, we obtain
\[
\xymatrix{
\Sigma^{(I)}\ar[rr]\ar@{=}[d]&&\Sigma^{(J)}\ar[rr]\ar@{=}[d]&&M\ot_R\Sigma\ar[rr]\ar[d]&&0\\
F(R)^{(I)}\ar[rr]&&F(R)^{(J)}\ar[rr]&&F(M)\ar[rr]&&0
}
\]
By the universal property of the cokernels $(M\ot_R\Sigma,g\ot_R\Sigma)$ and $(FM,Fg)$ we obtain that $F(M)\cong M\ot_R\Sigma$, this isomorphism is easily verified to be natural.

$\ul{(v)}$ follows from the fact that $\Sigma\in{_R\Mm^\cc}$, by construction.

$\ul{(vi)}$ follows by direct calculation.
\end{proof}

The following Theorem generalizes \cite[3.1]{Wis:galcom}.

\begin{theorem}
Let $R$ be a firm ring, $\cc$ an $A$-coring and $\Sigma\in{_R\Mm^\cc}$. Then the following assertions are equivalent.
\begin{enumerate}[(i)]
\item $\Sigma$ is an $R$-$\cc$ comonadic-Galois comodule;
\item For every $(\cc,A)$-injective $N\in\Mm^\cc$, the evaluation map
$$\ev_N:\Hom^\cc(\Sigma,N)\ot_R\Sigma\to N,\ \ev_N(f\ot_Ru)=f(u);$$
is an isomorphism
\end{enumerate}
\end{theorem}

\begin{proof}
$\ul{(i)\Rightarrow (ii)}$.
Recall that a right $\cc$-comodule $N$ is $(\cc,A)$-injective if and only if the comultiplication $\rho_N:N\to N\ot_A\cc$ has a left inverse $\gamma_N$ in $\Mm^\cc$.
For all $L\in \Mm^\cc$, consider the following diagram
\begin{equation}\eqlabel{equaliser}
\xymatrix{
\Hom^\cc(L,N)\ar[r]^-{i} & \Hom_A(L,N)\ar@<.5ex>[r]^-{j_1} \ar@<-.5ex>[r]_-{j_2} &\Hom_A(L,N\ot_A\cc).
}
\end{equation}
where the maps $j_1$ and $j_2$ are defined as follows
$$j_1(f)(l)=f(l)_{[0]}\ot_A f(l)_{[1]},\qquad j_2(f)(l)=f(l_{[0]})\ot_A l_{[1]}.$$
Define as well maps $\alpha:\ \Hom_A(L,N)\to \Hom^\cc(L,N)$
and $\beta:\ \Hom_A(L,N\ot_A\cc)\to \Hom_A(L,N)$ which are given by the formulas
$$\alpha(f)(l)=\gamma_N(f(l_{[0]})\ot_A l_{[1]}),\qquad
\beta(g)=\gamma_N\circ g.$$
Now it is easy to check that $j_1\circ i=j_2\circ i$, $\alpha\circ i=\Hom_A(L,N)$, $\beta\circ j_1=\Hom^\cc(L,N)$ and $\beta\circ j_2=i\circ \alpha$, this means that \equref{equaliser} is a contractable equaliser. Consequently, If we apply the functor $-\ot_R\Sigma$ to \equref{equaliser}, we obtain a contractable equaliser in $\Mm^\cc$ (see e.g. \cite[Proposition 3.3.2]{BarrWells:ttt} or \cite[VI.6]{McLane}). By the coassociativity condition, we can associate to $L$ a second equaliser in $\Mm^\cc$,
\[
\xymatrix{
L \ar[r]^{\rho_L} & L\ot_A\cc \ar@<.5ex>[rr]^-{\rho_L\ot_A\cc} \ar@<-.5ex>[rr]_-{L\ot_A\Delta} && L\ot_A\cc\ot_A\cc 
}
\]
Now take $L=\Sigma$, than we can compare both equalisers by the cannonical map,
\[
\xymatrix{
\Hom^\cc(\Sigma,N)\ot_A\Sigma \ar[d]_{\ev_N} \ar[r]^-{i} & \Hom_A(\Sigma,N)\ot_A\Sigma \ar[d]^{\can_N} \ar@<.5ex>[r]^-{j_1} \ar@<-.5ex>[r]_-{j_2} &\Hom_A(\Sigma,N\ot_A\cc)\ot_A\Sigma  \ar[d]^{\can_{N\ot_A\cc}}\\
\Sigma \ar[r]^{\rho_\Sigma} & \Sigma\ot_A\cc \ar@<.5ex>[r]^-{\rho_\Sigma\ot_A\cc} \ar@<-.5ex>[r]_-{\Sigma\ot_A\Delta} & \Sigma\ot_A\cc\ot_A\cc 
}
\]
Since $\can_N$ and $\can_{N\ot_A\cc}$ are isomorphisms, we find that $\ev_N$ is an isomorphism by the universal property of the equalisers.

$\ul{(ii)\Rightarrow (i)}$. For any right $A$-module $M$, the right $\cc$-comodule $M\ot_A\cc$ is $(A,\cc)$-injective. Consider the isomorphism $$\varpi_M=(-\ot_A\cc)\circ\rho_\Sigma:\Hom_A(\Sigma,M)\cong\Hom^\cc(\Sigma,M\ot_A\cc)$$
one can easily check that $\can_M=\ev_{M\ot_A\cc}\circ(\varpi_M\ot_A\Sigma)$, so we conclude that $\can_M$ is an isomorphism for all $M\in\Mm_A$ and $\Sigma$ is $R\hbox{-}\cc$ comonadic-Galois.
\end{proof}

\begin{theorem}\thlabel{equiv}
Let $R$ be a firm ring and $\cc$ an $A$-coring. Suppose that we have a pair of adjoint functors $(F,G)$ as follows
\[
\xymatrix{
\Mm_R \ar@<0.5ex>[rr]^F && \Mm^\cc \ar@<0.5ex>[ll]^G
}
\]
Put $\Sigma=F(R)$.
Then the following statements hold
\begin{enumerate}[(a)]
\item If $G$ has a right adjoint then
\begin{enumerate}[(i)]
\item we have a comatrix coring context $(R,A,\Sigma,\Sigma^\dagger,\sigma,\tau)$, with $\Sigma^\dagger=\Sigma^*\ot_RR$;
\item $\Sigma$ is $R$-firmly projective (as a right $\cc$-comodule);
\item if $R$ is flat as a left $R$-module then $G\cong-\otimes^\cc(\Sigma^*\ot_RR)$.
\end{enumerate}
\item If $(F,G)$ establish an equivalence of categories, then
\begin{enumerate}[(i)]
\item $\Sigma$ is an $R$-$\cc$ Galois comodule, this implies in particular $\cc\cong \Sigma^\dagger\ot_R\Sigma$ and $\*c\cong {_R\End}(\Sigma)$;
\item $\Sigma$ is an $R$-$\cc$ comonadic-Galois comodule, this implies in particular
$\cc\cong \Sigma^*\ot_R\Sigma$ and $\cc^*\cong\End_R(\Sigma^\dagger)$;
\item $\Sigma$ is a generator in $\Mm^\cc$; 
\item if $R$ is projective as a right $R$-module then $\Sigma$ is a projective generator in $\Mm^\cc$;
\item if $\cc$ is flat as left $A$-module then $\Sigma$ is faithfully flat as a left $R$-module.
\end{enumerate}
\end{enumerate}
%
Conversely, if $\cc$ is flat as a left $A$-module then $(b)(ii)$ and $(b)(vi)$ imply that\\
$(-\otimes_R\Sigma,\Hom^\cc(\Sigma,-)\ot_RR)$ is a pair of inverse equivalences between $\Mm_R$ and $\Mm^\cc$.
\end{theorem}

\begin{proof}
$\ul{(a)(i)\& (ii)}$. Consider the diagram of adjoint functors
\[
\xymatrix{
\Mm_R \ar@<0.5ex>[rrrr]^-{F} \ar@<0.5ex>[rrrrdd]^-{F'} &&&& \Mm^\cc \ar@<0.5ex>[llll]^-{G} \ar@<0.5ex>[dd]^-{F^\cc}\\
\\
&&&& \Mm_A \ar@<0.5ex>[uu]^-{G^\cc} \ar@<0.5ex>[lllluu]^-{G'}
}
\]
Where $G^\cc=-\ot_A\cc$ and $F^\cc$ is the forgetful functor. By \thref{Eilenberg-Wattscoring} we know that $F\cong -\ot_R\Sigma$ and $G\cong \Hom^\cc(\Sigma,-)\ot_RR$.
Then we can compute $F'=F^\cc\circ F=-\ot_A\Sigma$ and $G'=G\circ G^\cc=\Hom^\cc(\Sigma,-\ot_A\cc)\ot_RR\cong\Hom_A(\Sigma,-)\ot_RR$, and $(F',G')$ is a pair of adjoint functors. Denote the right adjoint of $G$ by $H$. Since $G^\cc$ has also a right adjoint, namely $H^\cc=\Hom^\cc(\cc,-)$, we can conclude that $G'$ has a right adjoint given by $H'=H^\cc\circ H$. By \thref{firmadjoint} we obtain immediately both statements, since $\Sigma\in{_R\Mm^\cc}$.

$\ul{(iii)}$ The statement follows now from part $(ii)$ and \cite[Theorem 4.4]{GTV}.


$\ul{(b)(i)}$. By part (a) we already know that $\Sigma$ is $R$-firmly projective as right $A$-module. Since $G$ is fully faithful we can obtain $\cc\cong FG(\cc)=\Hom^\cc(\Sigma,\cc)\ot_RR\ot_R\Sigma\cong \Sigma^\dagger\ot_R\Sigma\cong\Sigma^*\ot_R\Sigma$.

Let us know prove that $\*c\cong{_R\End}(\Sigma)$. Take $\varphi\in\*c$, then we define $\alpha:\*c\to\End(_R\Sigma)$ by 
$$\alpha(\varphi)(u)=e_r\varphi(f_r\otimes_Ru^r).$$
Conversely, we have a map $\beta:\End(_R\Sigma)\to \*c$  :
$$\beta(\psi)(f\otimes_Ru)=f(\psi(u)),$$
with $\psi\in\End(_R\Sigma)$.
We check that $\alpha$ and $\beta$ are each other inverses.
\[
\begin{array}{rcl}
\beta\circ\alpha(\varphi)(f\otimes_Ru)&=&f(e_r\varphi(f_r\otimes_Ru^r))=f(e_r)\varphi(f_r\otimes_Ru^r)\\
&=&\varphi(f(e_r)f_r\otimes_Ru^r)=\varphi(f\cdot r\otimes_Ru^r)\\
&=&\varphi(f\otimes_Rr\cdot u^r)=\varphi(f\otimes_Ru)\\
{}\\
\alpha\circ\beta(\psi)(u)&=&e_rf_r(\psi(u^r))=r\cdot\psi(u^r)\\
&=&\psi(r\cdot u^r)=\psi(u)\\
\end{array}
\]
Finally, let us check that $\alpha$ is a ring morphism. 
\[
\begin{array}{rcl}
\alpha(\phi *\varphi)(u)&=&e_r(\phi*\varphi)(f_r\otimes_Ru^r)\\
&=&e_r\varphi(f_r\otimes_Re_s\phi(f_s\otimes_Ru^{rs}))\\
{}\\
\alpha(\varphi)\circ\alpha(\phi)(u)&=&\alpha(\varphi)(e_r\phi(f_r\otimes_Ru^r))\\
&=&e_s\varphi(f_s\otimes_Re_r^s\phi(f_r\otimes_Ru^r))
\end{array}
\]
Both are equal since $e_r\otimes_Af_r\otimes_Re_s\otimes_Af_s\otimes_Ru^{rs}=e_s\otimes_Af_s\otimes_Re_r^s\otimes_Af_r\otimes_Ru^r$. 

$\ul{(ii)}$. Take any $M\in\Mm_A$, since $G$ is fully faithful we find,
$$M\otimes_A\cc\cong FG(M\otimes_A\cc)\cong\Hom^\cc(\Sigma,M\otimes_A\cc)\otimes_R\Sigma
\cong\Hom_A(\Sigma,M)\otimes_R\Sigma.$$
So $\Sigma$ is $R$-$\cc$ comonadic Galois.
Since $-\ot_R\Sigma$ and $\Hom^\cc(\Sigma,-)\ot_RR$ make up a pair of adjoint functors, for any $M\in\Mm_R$ and $N\in\Mm^\cc$, we obtain an isomorphism
$$\Hom^\cc(M\ot_R\Sigma,N)\cong\Hom_R(M,\Hom^\cc(\Sigma,N)\ot_RR).$$
When we apply this to the situation $\Sigma^\dagger\in\Mm_R$ and $\cc\in\Mm^\cc$, we find
$$\cc^*\cong\Hom^\cc(\Sigma^\dagger\ot_R\Sigma,\cc)\cong
\Hom_R(\Sigma^\dagger,\Hom^\cc(\Sigma,\cc)\ot_RR)\cong\End_R(\Sigma^\dagger).$$

$\ul{(iii)}$. This statement follows from classical arguments, but we give the proof for the sake of completeness. Take any $N\in\Mm^\cc$. Since we know that $R$ is a generator in $\Mm_R$, there exists an exact row (epimorphism) in $\Mm_R$ of the following form
\[\xymatrix{
R^{(I)}\ar[rr]&& G(N) \ar[rr] && 0
}\]
Since $(F,G)$ is an adjoint pair, $F$ is preserves epimorphisms and direct sums, so we obtain the following epimorphism in $\Mm^\cc$
\[\xymatrix{
F(R^{(I)})=\Sigma^{(I)}\ar[rr]&& FG(N)=N \ar[rr] && 0
}\]
where we used that $G$ is a full and faithful functor. It follows that $\Sigma$ is generator in $\Mm^\cc$. 

$\ul{(iv)}$ To prove that $\Sigma$ is projective, one can proceed in a similar way as in part $(iv)$. Suppose the following diagram with exact row (epimorphism) is given in $\Mm^\cc$
\[
\xymatrix{
&&\Sigma \ar[d]^{g}\\
N\ar[rr]^-{f} && N' \ar[rr] && 0
}
\]
After applying the functor $G$, which preserves epimorphisms, we can complete the diagram with a morphism $h$, since $R$ is projective as a right $R$-module.
\[
\xymatrix{
&&G(\Sigma)=R \ar[d]^{G(g)} \ar[dll]_{h}\\
G(N)\ar[rr]^-{G(f)} && G(N') \ar[rr] && 0
}
\]
Apply now the functor $F$, then we obtain 
\[
\xymatrix{
&&F(R)=\Sigma \ar[d]^{FG(g)=g} \ar[dll]_-{F(h)}\\
FG(N)=N\ar[rr]^-{FG(f)=f} && FG(N')=N' \ar[rr] && 0
}
\]
This shows that $\Sigma$ is a projective right $\cc$-comodule.

The remaining statements follow from the stucture theorem for $R$-$\cc$ Galois comodules, see \cite[Theorem 4.15]{GTV}.
\end{proof}

\begin{corollary}\colabel{locunits}
Let $R$ be a ring with local units and $\cc$ an $A$-coring. Suppose we have an equivalence of categories
\[
\xymatrix{
\Mm_R \ar@<0.5ex>[rr]^F && \Mm^\cc \ar@<0.5ex>[ll]^G
}
\]
Put $\Sigma=F(R)$. Then the following statements hold
\begin{enumerate}[(i)]
\item $\Sigma$ is a locally projective right $A$-module (in the sense of Zimmermann-Huisgen \cite{ZH});
\item $F\cong-\otimes_R\Sigma$ and $G\cong-\otimes^\cc\Sigma^\dagger$, with $\Sigma^\dagger=\Sigma^*\ot_RR$;
\item We have a comatrix coring context $(R,A,\Sigma,\Sigma^\dagger,\sigma,\tau)$ and $\cc$ is isomorphic to the comatrix coring $\Sigma^*\otimes_R\Sigma$;
\end{enumerate}
If $\cc$ is flat as a left $A$-module, then the following statement holds as well
\begin{enumerate}
\item[(iv)] $\Sigma$ is faithfully flat as a left $R$-module;
\end{enumerate}
Conversely, $(i)$, $(iii)$ and $(iv)$ imply that\\
$\cc$ is flat as a left $A$-module and $(-\otimes_R\Sigma,\Hom^\cc(\Sigma,-)\ot_RR)$ is a pair of inverse equivalences.
\end{corollary}

\begin{proof}
Since $R$ is a ring local units, $R$ is locally projective and consequently flat as a left $R$-module.
Statements $(ii)$, $(iii)$ and $(iv)$ follow directly from \thref{equiv}. Also by \thref{equiv}, we know that $\Sigma$ is a firm module over $R$, with action induced by a ring morphism $R\to \Sigma\ot_A\Sigma^\dagger$, where $(\Sigma,\Sigma^\dagger,\tau)$ constitute a dual pair. Since $R$ has local units, this implies $Z^\dagger=\Sigma\ot_A\Sigma^\dagger$ is also a ring with local units and $\Sigma$ is a firm left $Z^\dagger$-module. By \cite[Theorem 3.4]{V} we obtain that $\Sigma$ is locally projective as a right $A$-module.

For the last statement, we only have to prove that $\cc$ is flat as a left $A$-module. As in the previous part, we find that the firm right $R$-module $\Sigma^\dagger$ is locally projective as a left $A$-module, and consequently $\Sigma^\dagger$ is flat as a left $A$-module. Since we also suppose $\Sigma$ to be flat as a left $R$-module, we find that $\cc=\Sigma^\dagger\ot_R\Sigma\cong\Sigma^*\ot_R\Sigma$ is a flat left $A$-module.
\end{proof}

In the next corollary we use the notation and terminology of split direct systems that was introduced in \cite{CDV:colimit}.

\begin{corollary}\colabel{idlocunits}
Let $R$ be a ring with idempotent local units and $\cc$ an $A$-coring. Suppose we have an equivalence of categories
\[
\xymatrix{
\Mm_R \ar@<0.5ex>[rr]^F && \Mm^\cc \ar@<0.5ex>[ll]^G
}
\]
Put $\Sigma=F(R)$. Then the following statements hold
\begin{enumerate}[(i)]
\item $\Sigma\cong\colim \ul{P}$ for a split directed system $\ul{P}^s:\Zz\to {\Mm^\cc_{\rm fgp}}^s$
\item $F\cong-\otimes_R\Sigma$ and $G\cong-\otimes^\cc\Sigma^\dagger$, with $\Sigma^\dagger=\colim \ul{P^*}$;
\item We have a comatrix coring context $(R,A,\Sigma,\Sigma^\dagger,\sigma,\tau)$ and $\cc$ is isomorphic to the comatrix coring $\Sigma^*\otimes_R\Sigma$;
\item $R\cong \colim \ul{T}$ for a direct system $\ul{T}:\Zz\to \Ff_k$, where $T_i=\End^\cc(P_i)$ for $i\in \Zz$.
\end{enumerate}
If $\cc$ is flat as a left $A$-module, then the following statement holds as well
\begin{enumerate}
\item[(v)] $\Sigma$ is faithfully flat as a left $R$-module;
\end{enumerate}
Conversely, $(i)$, $(iii)$ and $(v)$ imply that\\
$\cc$ is flat as a left $A$-module and $(-\otimes_R\Sigma,\Hom^\cc(\Sigma,-)\ot_RR)$ is a pair of inverse equivalences.
\end{corollary}

\begin{proof}
We can argue as in the proof of \coref{locunits}, to obtain that $Z^\dagger=\Sigma\ot_A\Sigma^\dagger$ is a ring with idempotent local unit and $\Sigma$ and $\Sigma^\dagger$ are firm left and right $Z^\dagger$-modules.
We can apply \cite[Lemma 4.3, Lemma 4.11]{CDV:colimit}, to obtain a characterisation of $\Sigma$ and $\Sigma^\dagger$ in terms of a colimit of a split direct system. It follows from the proofs of \cite[Lemma 4.3, Lemma 4.11]{CDV:colimit} that the individual modules $P_i$ of the split direct system are obtained as $e_i\Sigma$, where $e_i$ is an idempotent local unit of $R$. Since $\Sigma\in{_R\Mm^\cc}$, the action of $e_i$ is right $\cc$-colinear and thus $P_i$ is a right $\cc$-comodule.

Since modules that are discribed as a colimit of a split direct system are in particular locally projective in the sense of Zimmermann-Huisgen (see \cite{V}), all other statements follow now from \thref{equiv} and \coref{locunits}.
\end{proof}

\begin{corollary}\colabel{equivunit}
Let $T$ be a ring with unit and $\cc$ an $A$-coring. Suppose we have an equivalence of categories
\[
\xymatrix{
\Mm_T \ar@<0.5ex>[rr]^F && \Mm^\cc \ar@<0.5ex>[ll]^G
}
\]
Put $\Sigma=F(T)$. Then the following statements hold
\begin{enumerate}[(i)]
\item $\Sigma$ is a finitely generated and projective right $A$-module;
\item $F\cong-\otimes_T\Sigma$ and $G\cong-\otimes^\cc\Sigma^*$;
\item We have a comatrix coring context $(T,A,\Sigma,\Sigma^*,\sigma,\tau)$ and $\cc$ is isomorphic to the comatrix coring $\Sigma^*\otimes_T\Sigma$;
\item the map $\jmath : T\to \End^\cc(\Sigma)$ is an isomorphism.
\end{enumerate}
If $\cc$ is flat as a left $A$-module, then the following statement holds as well
\begin{enumerate}
\item[(v)] $\Sigma$ is a finitely generated and projective generator in $\Mm^\cc$; 
\item[(vi)] $\Sigma$ is faithfully flat as a left $T$-module;
\end{enumerate}
Conversely, $(i)$, $(iii)$ and $(vi)$ imply that\\
$\cc$ is flat as a left $A$-module and $(-\otimes_T\Sigma,\Hom^\cc(\Sigma,-))$ is a pair of inverse equivalences with $T=\End^\cc(\Sigma)$.
\end{corollary}

\begin{proof}
By \thref{equiv} and the above corollaries, we only have to prove the generating part of statement $(v)$. This proof can be found in \cite[Theorem 3.2]{EGT:comcor}.
Recall that an object $X$ in an abelian category $\Cc$ is called finitely generated if and only if for any directed family of subobjects $\{X_i\}_{i\in I}$ of $X$ satisfying $X=\sum_{i\in I} X_i$, there exists an $i_0\in I$ such that $X=X_{i_0}$. If $X$ is a generator for $\Cc$, this condition is known to be equivalent with the fact that $\Hom_\Cc(X,-)$ preserves coproducts. Since $\cc$ is flat as left $A$-module, we know that $\Mm^\cc$ is a Grothendieck category and $\Hom^\cc(\Sigma,-)=G$, which clearly preserves coproducts.
\end{proof}

\section{Comonadic-Galois versus Galois comodules}\selabel{comparision}

\begin{proposition}\prlabel{connectiongal}
Let $\cc$ be an $A$-coring, $R$ a firm ring and $\Sigma\in{_R\Mm^\cc}$.
If $\Sigma$ is an $R\hbox{-}\cc$ Galois comodule, then $S\ot_S\Sigma\in{_S\Mm^\cc}$ is an $S$-$\cc$ comonadic-Galois comodule, for any firm ring $S$ such that $R\subset S\subset T=\End^\cc(\Sigma)$.
\end{proposition}

\begin{proof}
Remark that $\Sigma$ is not always a firm left $S$-module, even if it is firm as a left $R$-module, for this reason we have to consider $S\ot_S\Sigma$.
Suppose $\Sigma$ is a $R\hbox{-}\cc$ Galois comodule. 
We construct an inverse map for the natural morphism \equref{nattfgalois}. Define
\[
\nu_M:M\otimes_A\cc \to \Hom_A(\Sigma,M)\otimes_SS\ot_S\Sigma,
\]
as the composite of the following morphisms
\[
\hspace{-1cm}
\xymatrix{
M\otimes_A\cc \ar[rr]^-{M\otimes_A\can^{-1}} &&M\otimes_A\Sigma^*\otimes_R\Sigma \ar[r] &\Hom_A(\Sigma,M)\otimes_RR\ot_R\Sigma \ar[r]& \Hom_A(\Sigma,M)\otimes_SS\ot_S\Sigma,
}
\]
We use the following notation
${\can^{-1}}:\cc\to\Sigma^*\otimes_R\Sigma,~\can^{-1}(c)=f_c\otimes_Au_c$, then we check
\begin{equation}\eqlabel{vgl8}
\nu_M\circ\can_M(\phi\otimes_Ss\ot_S u)=\nu_M(\phi s(u_{[0]})\otimes_Au_{[1]})=
\phi s(u_{[0]})f_{u_{[1]}}\otimes_S r\ot_S (u_{u_{[1]}})^r
\end{equation}
Since $\Sigma\in{_R\Mm^\cc}$, the following diagram commutes
\begin{equation*}
\xymatrix{
\Sigma \ar[rr]^{\rho_{\Sigma,\cc}} \ar[d]_{\rho_{\Sigma,\dd}} && \Sigma\ot_A\cc \\
\Sigma\ot_A\Sigma^*\ot_R\Sigma \ar[rru]_{\Sigma\ot_A\can}
}
\end{equation*}
i.e. $\rho_{\Sigma,\cc}(u)=u_{[0]}\ot_Au_{[1]}=e_rf_r(u^r_{[0]})\ot_Au^r_{[1]}$. Applying $\Sigma\ot_A\can^{-1}$ to, we obtain $u_{[0]} \ot_A f_{u_{[1]}}\otimes_R u_{u_{[1]}}=e_r\ot_Af_r\ot_Ru^r$. If we apply now the firmness property of $\Sigma$ on the last factor in the tensor product, we find
$u_{[0]} \ot_A f_{u_{[1]}}\otimes_R r\ot_R (u_{u_{[1]}})^r =e_r\ot_Af_r\ot_Rr'\ot u^{rr'}$
With this equality we can rewrite the last formula of \equref{vgl8} in the following way
\begin{eqnarray*}
\phi s(u_{[0]})f_{u_{[1]}}\otimes_Sr\ot_S (u_{u_{[1]}})^r
&=&\phi s(e_r)f_r\otimes_Sr'\ot_Su^{rr'}\\
&=&\phi s\cdot r\otimes_Sr'\ot_S u^{rr'}=\phi s\otimes_Sr\ot r'u^{rr'}\\
&=&\phi\otimes_Ss\ot_S u
\end{eqnarray*}

For the converse, first remark that 
$\varepsilon_\cc=\ev\circ\can_A^{-1}$, where $\ev$ is the evaluation map
$$\ev:\Sigma^*\otimes_S\Sigma\to A,\qquad\ev(\varphi\otimes_Au)=\varphi(u).$$
This follows from the straightforward computation
$$\varepsilon_\cc\circ\can_A(\varphi\ot_Su)=\varepsilon(\varphi(u_{[0]})u_{[1]})=\varphi(u)=
\ev(\varphi\ot_Su),$$
and from $\can\circ\can^{-1}=\cc$. Apply this new identity in the penultimate equality of the next computation
\begin{eqnarray*}
\can_M\circ\nu_M(m\otimes_A c)&=&\can_M(mf_c\otimes_S r\ot_S(u_c)^r)=
mf_c r((u_c)^r_{[0]})\otimes_A(u_c)^r_{[1]}\\
&=&mf_c ((r(u_c)^r)_{[0]})\otimes_A (r(u_c)^r)_{[1]}
= mf_c ((u_c)_{[0]})\otimes_A (u_c)_{[1]}\\
&=&m\otimes_Af_c((u_c)_{[0]})(u_c)_{[1]}=m\otimes_Af_{c_{(1)}}(u_{c_{(1)}})c_{(2)}\\
&=&m\otimes_A\varepsilon_\cc(c_{(1)})c_{(2)}=m\otimes_Ac,
\end{eqnarray*}
where we used in the third equality that $r$ is colinear and in sixth equality that $\can$ a morphism of right comodules.
\end{proof}

\begin{theorem}\thlabel{comparmain}
Let $\cc$ be an $A$-coring, $R$ a firm ring and $\Sigma\in{_R\Mm^\cc}$. If $\Sigma$ is $R$-firmly projective, then $\Sigma$ is an $R\hbox{-}\cc$ Galois comodule, if and only if $\Sigma$ is an $R$-$\cc$ comonadic-Galois comodule.
\end{theorem}

\begin{proof}
If $\Sigma$ is an $R$-$\cc$ Galois comodule, then it is an $R$-$\cc$ comonadic-Galois comodule by \prref{connectiongal}.  
To prove the converse, remember that since we know that $\Sigma$ is $R$-firmly projective, we can construct the comatrix coring $\Sigma^*\ot_R\Sigma$. Moreover, $\Sigma$ is comonadic-Galois, so in particular $\can_A : \Sigma^*\ot_R\Sigma\to \cc$ is an isomorphism (of right $\cc$-comodules). This map is exactly the canonical coring morphism $\can$ \equref{can}, so $\Sigma$ is also an $R$-$\cc$ Galois comodule.
\end{proof}

\begin{theorem}\thlabel{comparision3}
Let $\cc$ be an $A$-coring, $R$ a firm ring and $\Sigma\in{_R\Mm^\cc}$. We denote $T=\End^\cc(\Sigma)$
\begin{enumerate}[(i)]
\item
If $\Sigma$ is an $R$-$\cc$ comonadic Galois module, then $\Sigma$ is also an $T$-$\cc$ comonadic Galois comodule (i.e. $\Sigma$ is a Galois comodule in the sense of Wisbauer). 
\item 
If $R$ is a left ideal in $T$ and $\Sigma$ is an $T$-$\cc$ comonadic Galois module then $\Sigma$ is also an $R$-$\cc$ comonadic Galois comodule (and thus both properties are equivalent). 
\end{enumerate}
\end{theorem}

\begin{proof}
$\ul{(i)}$ Follows immediately form the commutativity of the following diagram and the surjectivity of $\pi$,
\[
\xymatrix{
\Hom_A(\Sigma,M)\ot_S\Sigma \ar[rr]^-{\can_{M,T}} && M\ot_A\cc \\
\Hom_A(\Sigma,M)\ot_R\Sigma \ar[rru]_{\can_{M,R}} \ar[u]^{\pi_M} 
}
\]
$\ul{(ii)}$ If $R$ is a left ideal in $T$, then we know by \cite[Lemma 4.11]{GTV} that 
$$\pi_M: \Hom_A(\Sigma,M)\ot_R\Sigma\to \Hom_A(\Sigma,M)\ot_T\Sigma$$ 
is an isomorphism for all $M\in\Mm_A$. Consequently $\can_{M,T}$ is an isomorphism if and only if $\can_{M,R}$ is an isomorphism.
\end{proof}

\begin{remark}
That the converse of statement (i) of \thref{comparision3} does not hold in general, follows from the following example. Let $\cc=A$ be the trivial $A$-coring, and take $\Sigma=A$. Then $\End_A(A)=A$ and for all $M\in\Mm_A$, the cannonical $\can_{M,A}$ is the trivial isomorphism $\can_{M,A}: \Hom_A(A,M)\ot_A A\to M\ot_A A$. Take any nontrival ringmorphism $R\to A$, such that $R$ is no ideal in $A$. Then we can compute the maps $\pi_M$ and $\can_{M,R}$ as
$$\Hom_A(A,M)\ot_BA\cong M\ot_BA \to M\ot_A A\cong M$$
Which are clearly not isomorphisms in general.

A sufficient condition to obtain that $R$ is an ideal in $T$, and thus that the Galois-property in the sense of Wisbauer \cite{Wis:galcom} is equivalent with $R$-$\cc$ comonadic Galois is that the functor $G=\Hom^\cc(\Sigma,-)\ot_RR$ \equref{expressionG} is a full and faithful functor (see \cite[Theorem 4.15, $(iv)\Rightarrow(v)$]{GTV} or \cite[Lema 2.4]{GT:comonad}).
\end{remark}

%
%

Take any $\Sigma\in{_S\Mm^\cc}$, for any firm ring $S\subset T=\End^\cc(\Sigma)$. As before, let $Z=\Sigma\otimes_A\Sigma^*$ be the elementary (possibly non-unital) $S$-ring associated to the dual pair $(\Sigma,\Sigma^*,\ev)$.
If moreover $\Sigma$ is a $S$-$\cc$ comonadic-Galois comodule, then we can also construct comultiplications on $Z$ and $\Sigma$ as follows.
\begin{eqnarray*}
\sd_{Z,\Sigma}&=&(\Sigma\ot_A\can_A^{-1})\circ\rho_{\Sigma,\cc}\\
\sd_Z&=&\sd_{Z,\Sigma}\ot_A\Sigma^*
\end{eqnarray*}
Remark that $\sd_{Z,\Sigma}$ is equal to the map that is obtained by composing the obvious map $\Sigma\to \End_A(\Sigma)\ot_S\Sigma,u\mapsto \End_A(\Sigma)\ot_Su$, with the isomorphism 
$$(\Sigma\ot_A\can^{-1}_A)\circ\can_\Sigma:\End_A(\Sigma)\otimes_S\Sigma \to \Sigma\ot_A\Sigma^*\otimes_S\Sigma.$$

Consider the following diagram
\[\hspace{-1cm}
\xymatrix{
\Sigma \ar[rr]^-{\rho_\Sigma} \ar[dd]_-{\rho_\Sigma}
&& \Sigma\ot_A\cc \ar[rr]^-{\Sigma\ot_A\can^{-1}} \ar[dd]_-{\rho_\Sigma\ot_A\cc}
&& \Sigma\ot_A\Sigma^*\ot_S\Sigma \ar[dd]_-{\rho_\Sigma\ot_A\Sigma^*\ot_S\Sigma}\\
&(1)
 && (2)\\
\Sigma\ot_A\cc \ar[rr]^-{\Sigma\ot_A\Delta} \ar[dd]_-{\Sigma\ot_A\can^{-1}}
&& \Sigma\ot_A\cc\ot_A\cc \ar[rr]^-{\Sigma\ot_A\cc\ot_A\can^{-1}} 
\ar[dd]_-{\Sigma\ot_A\can^{-1}\ot_A\cc}
&& \Sigma\ot_A\cc\ot_A\Sigma^*\ot_S\Sigma \ar[dd]_-{\Sigma\ot_A\can^{-1}\ot_A\Sigma^*\ot_S\Sigma}\\
& (3) 
&& (4)\\
\Sigma\ot_A\Sigma^*\ot_S\Sigma \ar[rr]^-{\Sigma\ot_A\Sigma^*\ot_S\rho_\Sigma}
&& \Sigma\ot_A\Sigma^*\ot_S\Sigma\ot_A\cc \ar[rr]^-{\Sigma\ot_A\Sigma^*\ot_S\Sigma\ot_A\can^{-1}}
&& \Sigma\ot_A\Sigma^*\ot_S\Sigma\ot_A\Sigma^*\ot_S\Sigma
}
\]
The commutativity of this diagram can be checked as follows. The quadrangle (1) commutes by coassociatityty of $\Sigma$, the commutativity of (2) and (3) follow by direct computation (use $\can$ in stead of $\can^{-1}$) and the commutativity of $(4)$ is trivial. We obtain that
\begin{eqnarray*}
(\sd_Z\ot_S\Sigma)\circ\sd_{Z,\Sigma}&=&
(Z\ot_S\sd_{Z,\Sigma})\circ\sd_{Z,\Sigma},
\end{eqnarray*}
after tensoring the above equality with $\Sigma^*$ we find as well
\begin{eqnarray*}
(\sd_Z\ot_SZ)\circ\sd_Z&=&
(Z\ot_S\sd_Z)\circ\sd_Z.
\end{eqnarray*}
So $Z$ is a (non-counital) $S$-coring and $\Sigma$ is a left $Z$-comodule. Moreover since $\varepsilon_\cc=\ev\circ\can^{-1}$, we have
\begin{eqnarray}\eqlabel{Sigmafirm1}
\mu_{Z,\Sigma}\circ\sd_{Z,\Sigma}=\Sigma;\\
\mu_Z\circ\sd_Z=Z.\nonumber
\end{eqnarray}
But we have more

\begin{proposition}\thlabel{compar2}
Let $\cc$ be an $A$-coring and $S$ any firm ring. Take $\Sigma\in{_S\Mm^\cc}$.
The following statements are equivalent
\begin{enumerate}
\item
The following canonical map is an isomorphism of $A\hbox{-}A$ bimodules
$$\can:\Sigma^*\ot_S\Sigma\to \cc,\qquad \can(f\ot_Su)=f(u_{[0]})u_{[1]};$$ (e.g. $\Sigma$ is a $S$-$\cc$ comonadic-Galois comodule) 
\item
we have comatrix coring context $(Z,A,\Sigma,\Sigma^\dagger,\eta,\epsilon)$ with $\Sigma^\dagger=\Sigma^*\otimes_SZ$ (i.e. $\Sigma$ is firmly projective) and 
$\cc$ is isomorphic to the associated comatrix coring.
\end{enumerate}
If there exists a firm ring $R$ and a ring morphism $\iota : R\to Z$, such that 
$\Sigma\in{_R\Mm^\cc}$ and $R$ is a left ideal in $S$, then
$\Sigma$ is an $R$-$\cc$ Galois comodule. 
\end{proposition}

\begin{proof}
Suppose first that $\can$ is an isomorphism. Let us prove that $\Sigma$ is firmly projective. Consider the maps $\mu_Z, \mu_{Z,\Sigma}, \sd_Z$ and $\sd_{Z,\Sigma}$ introduced above, where $Z$ is considered as a (nonunital) $S$-ring and a (non-counital) $S$-coring. Denote by $\bar{\mu}_Z:Z\ot_ZZ\to Z$ and $\bar{\sd}_Z:Z\to Z\ot_SZ\to Z\ot_ZZ$ the induced multiplication and comultiplication map, and similar for $\bar{\mu}_{Z,\Sigma}$ and $\bar{\sd}_{Z,\Sigma}$

By \equref{Sigmafirm1}, $\sd_{Z,\Sigma}$ is a right inverse of $\mu_{Z,\Sigma}$, and consequently also $\bar{\mu}_{Z,\Sigma}\circ\bar{\sd}_{Z,\Sigma}=\Sigma$. Let us show that $\bar{\sd}_{Z,\Sigma}$ is also a left inverse for $\bar{\mu}_{Z,\Sigma}$. Denote $\bar{\sd}_{Z,\Sigma}(u)=e_z\ot_Af_z\ot_Zu^z$. Since $\Sigma\in{_Z\Mm_A}$, we find
\begin{eqnarray*}
\bar{\sd}_{Z,\Sigma}\circ\bar{\mu}_{Z,\Sigma}(u\ot_A\varphi\ot_Zv)&=&
e_z\ot_Af_z\ot_Ru^z\varphi(v)= e_z\ot_Af_z\ot_Z(u^z\ot_A\varphi)\cdot v\\
&=&e_z\ot_Af_z\cdot(u^z\ot_A\varphi)\ot_Z v=e_z\ot_Af_z(u^z)\varphi\ot_Z v\\
&=&e_zf_z(u^z)\ot_A\varphi\ot_Z v=u\ot_A\varphi\ot_Z v
\end{eqnarray*}
So we find $\bar{\sd}_{Z,\Sigma}\circ\bar{\mu}_{Z,\Sigma}=\Sigma$, consequently $\bar{\sd}_Z\circ\bar{\mu}_Z=Z$, so $Z$ is a firm ring and $\Sigma$ is a firm left $Z$-module, i.e. $\Sigma$ is firmly projective. By \thref{firmadjoint} this implies that we have a comatrix coring context $(Z,A,\Sigma,\Sigma^\dagger,\eta,\epsilon)$ if we define
$$\epsilon:\Sigma^\dagger\ot_Z\Sigma\cong \Sigma^*\ot_S\Sigma \to A,\qquad \epsilon(\varphi\ot_Sz\ot_Zu)=\varphi(zu);$$
$$\eta=\bar{\sd}_Z:Z\to\Sigma\ot_Z\Sigma^\dagger=\Sigma\ot_A\Sigma^*\ot_SZ$$
One can easily check that $\can$ is a coring isomorphism between the associated comatrix coring and $\cc$.

The implication $(2) \Rightarrow (1)$ is trivial.

Finally, it follows from \cite[Lemma 4.11]{GTV} that $\cc\cong\Sigma^*\otimes_S\Sigma\cong\Sigma^*\otimes_R\Sigma$,
and $\Sigma$ is an $R\hbox{-}\cc$ Galois comodule.
\end{proof}

\section*{Acknowledgment}
The author would like to thank S. Caenepeel and G. B\"ohm for helpful comments on a previous version of this paper.

\end{document}